\input amstex
\documentstyle{ amsppt}
\magnification=\magstep1
\advance\hoffset by 0.1 truein
\hoffset=0in
\hsize=5in
\voffset=0in
\vsize=7in

\NoRunningHeads
\NoBlackBoxes

\topmatter
 \title On homaloidal polynomials \endtitle
 \author Andrea Bruno \endauthor
\endtopmatter
Let $\Bbb P^n$ be the projective space over a field of characteristic zero.
If $F$ is a homogeneous polynomial we say that $F$ is homaloidal
if the polar map $ \partial F$ defined by the partial derivatives of $F$ is a birational
selfmap of $\Bbb P^n$. Although the problem of determining
homaloidal polynomials has a classical flavour the theme only
recently was raised in an algebro-geometric context by Dolgachev
([Do]), following suggestions stemming from the theory of
prehomogeneous varieties: the relative invariants of
prehomogeneous spaces are in fact homaloidal polynomials
([KiSa],[EKP]). Dolgachev classifies homaloidal polynomials in $
\Bbb P^2$ (see also [Di]) and characterizes homaloidal polynomials in $\Bbb P^3$
which are products of linear forms as products of four independent
linear forms. Dolgachev also raises the question if it is true
that a non square free product of linear forms is homaloidal if
and only the product of its factors with multiplicity one is.
This question has been given a positive answer in a specific
case (see [KS]) and in full generality
([DiPa]) in a topological context.

We will give an algebraic proof of  the following: 
\proclaim {Theorem A} 
Suppose $F=
\prod_{i=0}^rL_i$ is a square free homaloidal polynomial which is
the product of linear forms. Then $r=n$ and the linear forms $L_0,\ldots,L_n$ are
independent, so that $\partial F$ is a composition of a
projectivity and of a standard Cremona transformation.
\endproclaim

\proclaim {Theorem B} Suppose $F=  \prod_{i=0}^rL_i^{n_i}$ is a
homogeneous polynomial with $L_0,\ldots,L_r$ linear forms.
Then $F$ is homaloidal if and only if the polynomial $F_{red}=
\prod_{i=0}^rL_i$ is homaloidal. In particular, if $F$ is
homaloidal, then $r=n$, the linear forms  $L_0,\ldots,L_n$ are independent and the
map  $\partial F$ is a composition of projectivities and of a
standard Cremona transformation.
\endproclaim

I would like to thank I. Dolgachev and F. Russo, for having brought
to my attention this problem, and M. Mella, for conversations on the subject.
This work has been written during and immediately after my stay in University of Catania,
during PRAGMATIC conference; I thank the organizers and S. Verra, who insisted for me to
publish it.
\vskip 0.5 cm
{\bf 1.Preliminaries}
We start with the following classical:\newline
{\bf Definition 1}
If $ F \in H^0(\Bbb P^n,\Cal O_{P^n}(d)) $ is a homogeneous
polynomial, the { \it polar map defined by \/} $F$ is the rational
map defined by $$ \partial F: \Bbb P^n - -\to \Bbb P^{n *} $$
$$ \partial F( p )= [ \frac {\partial F}{\partial X_{0}}(p), \ldots,  \frac {\partial F}{\partial X_{n}}(p)],$$
while the {\it polar system  defined by \/} $F$ is the linear system
$$| \partial F|:= |< \frac {\partial F}{\partial X_{0}}, \ldots,  \frac {\partial F}{\partial X_{n}}>| \subset
|\Cal O_{P^n}(d-1)|.$$ 
It follows from the definition that, if $ Z_F \subset \Bbb
P^n$ is the hypersurface defined by $F$, the base locus of the
polar map defined by $F$ is the singular locus of $Z_F$. Moreover,
since $Z_F$ is an hypersurface, $Z_F$ is not reduced if and only if
$F$ is not square free. In case $F$ is square free the polar map $\partial F$
is then free of base divisors and if $Z_F$ is smooth the polar map $\partial F$
is a morphism whose image is the dual
variety $Z_F^{*} \subset \Bbb P^{n *}$ of $Z_F$.
\vskip 0.2 cm
If $F$ is not square free, let us write $$ F =
\prod_{i=0}^rG_i^{n_i},$$
with $d= \sum_{i=0}^r n_i \text {deg} G_i$. Then
the base
components of the polar system $| \partial F|$ are given by the hypersurface defined
by the polynomial $$F^{'}= \prod_{i=0}^rG_i^{n_i-1}.$$
We will indicate by $F_{red}$ the polynomial $\frac {F}{F^{'}}$.
\vskip 0.2 cm
The polar system defined by $F$ is
naturally split in a fixed and in a moving part: 
{ \bf Definition 2} The { \it moving part \/} of the polar system
defined by a homogeneous polynomial $F$ is the linear system $ |
 M(\partial F)|$ obtained by removing all base components from the polar 
system $| \partial F|$.
\vskip 0.2 cm
In particular, we
have that $$ |M( \partial F)| \subset | \Cal O_{\Bbb P^n}(d-1- \text {deg} F^{'})|
.$$ Notice that the above definition makes perfectly
sense in case $F$ is square free, in which case $F^{'}$ is a
constant.

{\bf Definition 3}
A {\it homaloidal polynomial \/} of degree $d$ is a homogeneous polynomial
$ F \in H^0(\Bbb P^n, \Cal O_{P^n}(d))$ such that the moving
part $|M( \partial F)|$ of the polar system defined by $F$ induces a {\it birational map \/}.
\vskip 0.2 cm
Well known examples of homaloidal polynomials in $\Bbb P^n$ are those defining smooth quadrics. A remarkable result in [EKP] is the classification of homaloidal polynomials
of degree $d=3$: the irreducible ones define the secant varieties of the four Severi varieties and a classification seems at hand at least in degree $d=4$.
Probably the most known and important example of homaloidal polynomial is the
polynomial $F=X_0 \cdots X_n$, which has degree $d=n+1$ and whose
associated polar map is a {\it standard Cremona transformation \/}.
Another class of examples of homaloidal polynomials is given by the polynomials 
$$F(m_0,\ldots,m_n):= X_0^{m_0}\cdots X_n^{m_n},$$
with $m_i \geq 1$, for all $i=0, \ldots, n$,
of arbitrary large degree $d= \sum_{i=0}^n m_i$. The base
locus of $\partial F$ has a divisorial component defined by the
polynomial $F^{'}= \prod_{i=0}^nX_i^{m_i-1}.$ Once we remove it,
we simply compute that $$|M( \partial F)|=|<m_0\prod_{i=1}^nX_i,
\ldots,m_j\prod_{i \neq j }X_i, \ldots, m_n
\prod_{i=0}^{n-1}X_i>|,$$ so that $ \partial F$, induces the same
map as a composition of a projectivity (a homothety) and $ \partial F_{red}
 = \partial \prod_{i=0}^nX_i$, so that it is homaloidal.
\vskip 0.2 cm
To say that $F$ is homaloidal is equivalent to say that $\partial F$ is dominant and that there exists a resolution of singularities
$$
\CD
  \,    X       \\
\llap {$f$} \swarrow  \qquad \qquad \searrow \rlap {$g$} \\
\Bbb P^n --{\partial F}--\to    \Bbb P^{ n *},
\endCD
$$
such that if $Y$ is a general member of $|M( \partial F)|$ and if $\overline {Y}$ denotes its strict transform on $X$, we have $$( \overline {Y})^n=1,$$ because in fact  $ \overline {Y} \simeq g^{*} \Cal O_{\Bbb P^{n *}}(1) )$, where $\simeq$ denotes linear equivalence as Cartier divisors and $g$ is a birational morphism.
\vskip 0.2 cm
An important property of homaloidal polynomials is the following:
\proclaim { Proposition 4} If $F$ is a homaloidal polynomial,
$Z_F$ is not a cone. In particular if $F$ is the product of linear
forms $F= \prod_{i=0}^rL_i^{m_i}$, then  $ < L_0, \ldots, L_r > =
H^0(\Bbb P^n, \Cal O_{\Bbb P^n}(1))$, so that in particular $ r
\geq n$.
\endproclaim
\demo {Proof} If $Z_F$ is a cone the image of the polar map defined by $F$
is contained in the linear space dual to the vertex of the cone of $Z_F$.
 If $F$ is a product of
linear forms  $Z_F$ is a cone if and only if $<L_0,\ldots,L_r>
\neq H^0(\Bbb P^n, \Cal O_{\Bbb P^n}(1))$.\qed \enddemo
In fact (see [R]) even more is true: $Z_F$ is a cone if and only if the image of the polar map associated to $F$ lies in a hyperplane.
\vskip 0.5 cm

{\bf 2. Products of linear forms} In this section we will prove
Theorems A and B. We will always assume that $F$ is a homaloidal
polynomial which splits in the product of linear forms. Suppose
first that $F$ is square free so that $\partial F$ is a linear
system free of base components. We will fix a minimal resolution of
singularities:
$$
\CD
  \,    X       \\
\llap {$f$} \swarrow  \qquad \qquad \searrow \rlap {$g$} \\
\Bbb P^n --{\partial F}--\to    \Bbb P^{ n *}.
\endCD
$$
By definition, if $Y$ is a general member of the polar
system $|\partial F|$, and if $ \overline {Y}$ is its strict transform on $X$, 
we have $\overline {Y} \simeq g^{*} \Cal O_{\Bbb P^{n*}}(1).$

Let us write $$F=  \prod_{i=0}^rL_i.$$ 
Up to a projectivity we can assume that $L_0=X_0$. We will denote by $H_0$ the hyperplane 
defined by $X_0=0$. Let us define the polynomial $$G:= \frac {F}{X_0}= \prod _{i=1}^r L_i.$$
 With this choice, a basis of the polar system $| \partial F|$ defined by $F$ is given by:

$$ |\partial F|= | <G+ X_0 \frac {\partial G}{\partial X_0}, X_0 \frac 
{\partial G}{\partial X_1} , \ldots, X_0 \frac {\partial G}{\partial X_n}>|.$$
We first observe that, if $D \subset X$ is the strict transform of
the hyperplane $H_0$, looking at the equations for
$\partial F$, the map $g$ contracts $D$, because $\partial F$ contracts $H_0$ to its 
dual point $U_0=[1:0: \cdots:0]$.
\vskip 0.2 cm
Let us define $G_0 \in H^0(H_0, \Cal O_{H_0}(d-1))$ as the restriction of $G$ to $H_0$; 
notice that $G_0$ doesn't need to be reduced.
\proclaim {Lemma 5}
The irreducible divisor $D \subset X$ is the unique divisor contracting 
to the point $U_0 \in \Bbb P^{n*}$ and the map $g : X \to \Bbb P^{n*}$ 
factors through the blow up $h^{'} : Z \to \Bbb P^{n*}$ of $\Bbb P^{n*}$ at $U_0$.
\endproclaim
\demo {Proof}
Suppose that $W \neq D$ is a divisor in $X$ which is $g-$exceptional and such that 
$g(W)=g(D)=U_0$. By minimality of the resolution of the rational map $\partial F$, 
$W$ is not $ f-$exceptional, so that it corresponds to a hypersurface 
$f(W) \subset \Bbb P^n$ which is distinct from $H_0$. Let $J$ be an equation of 
$f(W)$. Looking at the equations of $ \partial F$, it follows that $J$ must divide 
$ \frac {\partial G}{ \partial X_i}$ for all $ i \geq 1$. If we restrict the system 
$|< \frac {\partial G}{ \partial X_1}, \ldots, \frac {\partial G}{ \partial X_n}>| 
$to the hyperplane $H_0$, by analyticity of polynomials, we obtain nothing but the 
polar system defined by $G_0$ on $H_0$.
This implies that $J$ restricted to $H_0$ is a base component of $ |\partial G_0|$ 
and this means that $$f(W) \cap H_0 \subset \cup \{ L_i | \, \text { there exists} \,
L_j \, \text { for which} \, L_j \in < H_0,L_i> \}.$$
But the polar map $\partial F$ contracts each one the hyperplanes defined by the 
linear forms $L_i$ to their dual points in $\Bbb P^{ n *}$, so that it cannot be 
$g(W)=g(H_0)$.

The irreducible divisor $D$ corresponds then to the extraction of a valuation centered at 
$U_0 \in \Bbb P^{n*}$ and we need to show that this valuation corresponds to the whole maximal ideal $m_{U_0}$.
We have already noticed en passant that the system $|< \frac {\partial G}{ \partial X_1}, 
\ldots, \frac {\partial G}{ \partial X_n}>| $ corresponds on $X$ to the system $ | g^{*} 
\Cal O_{\Bbb P^{n*}}(1)-D| $ and that it {\it is of codimension one \/} in $ | g^{*} \Cal 
O_{\Bbb P^{n*}}(1)| $.
Suppose that $g_{*} \Cal O_X(-D) = m^{'}$ with $\sqrt {m^{'}}= m_{U_0}$; since
$$g_{*} g^{*} \Cal O_{\Bbb P^{n*}}(1) \otimes \Cal O_X(-D)= m^{'} \otimes \Cal O_{\Bbb P^{n*}}(1),$$
we have that $H^0(\Bbb P^{n*}, m^{'} \otimes \Cal O_{\Bbb P^{n*}}(1))=n$ and hence that 
$$ m_{U_0}=m^{'}.$$ 
Hence $D$ is the strict transform of the exceptional divisor under the blow up
of $\Bbb P^{n*}$ at $U_0$, $h^{'}: Z \to \Bbb P^{n*}$ and the result follows\qed \enddemo

Consider now the following diagram of maps:

$$
\CD
  \,    X       \\
  \llap {$f$} \swarrow  \qquad \qquad \searrow \rlap {$g$} \\
\, \qquad H_0 \subset \Bbb P^n --{\partial F}--\to    \Bbb P^{ n *}-\pi - \to  P, \\
\endCD
$$
where $\pi$ is the projection from the point $U_0$ to
the hyperplane $P$. 
\vskip 0.2 cm 
We can use Lemma 5 in order to factorize the morphism $g$ through the blow up $Z$ of $ \Bbb P^{n *}$ at $U_0$. We have the following diagram:
$$
\CD
X \\
@VhVV \\
Z \\
\llap {$h^{'}$} \swarrow \qquad \qquad \searrow \rlap {$t$}\\
\Bbb P^{n *} - \pi - - \to P,
\endCD
$$
with $g=h^{'}h$.
Let us denote by $Y^{'}$ a general member of
the polar system $| \partial G|$  and recall that we denote by $Y$ a general element in  
$|\partial F|$ and by $D$ the strict
transform of $H_0$ in $X$.
\proclaim {Lemma 6}
With the above
notations:
\roster
\item the composition $th: X  \to P $ is a morphism and $(th)^{*} \Cal O_P(1) \simeq  
\overline {Y^{'}}  \simeq \overline { Y} -D$,
\item the restriction of the linear system $|(th)^{*} \Cal O_P(1)|$  to $D$ induces a morphism 
$m: D \to P$ which is a resolution of singularities of the polar map associated to $G_0$ on 
$H_0$, i.e. $|(\overline { Y} -D)|_D| = |M ( \partial G_0)|$, 
\item the polynomial $G_0$ is homaloidal in $H_0$,
\item $G_0$ is square free if and only there don't exist $L_i$ and $L_j$, with 
$i \neq j $, such that $X_0
\in <L_i,L_j>$
\endroster
\endproclaim

\demo {Proof }

It follows from the definitions that the morphism $th$ is set theoretically the same as the rational map $ \pi g$,
 so that in fact $th$ is a resolution of singularities of $\pi g$: they define the same morphism
 up to removing base components from the linear system defining $ \pi g$, which is in fact $D$.
 Now, the linear system defining $th$ is $|\overline {Y}-D|=|(th)^{*} \Cal O_P(1)|=|g^{*}\Cal O_{\Bbb P^{n*}}(1)-D|$; by the above argument, this is also the moving part of $|\partial G|$, so that in fact $Y^{'} \simeq \overline {Y}-D$.
It follows from the explicit equations that if we restrict to $H_0$ the map 
$ \pi g$, the resulting restriction is the rational map $ \partial G_0$, whose 
resolution of singularities is the map $th$ restricted to $D$.
We just notice here that the linear system inducing $m:= th|_D$ is  $|(\overline { Y} -D)|_D| = |M ( \partial G_0)|$. The map $\partial G_0$ 
is surjective because a composition of surjections.
In order to show that $G_0$ is homaloidal on $H_0$ it suffices to show that
$$D \cdot 
((th)^{*}\Cal O_P(1) -D)^{n-1}=1.$$ 
This follows from the fact that $D$ is the strict transform of the exceptional divisor
under the blow up $h^{'}: Z \to \Bbb P^{n*}$.
Finally, since the base components of $|\partial G_0|$ on $H_0$ are non reduced components of $Z_{G_0}$, and since $G$ is square free, the last part of the thesis is proved.\qed\enddemo 

{\bf Remark}
As a result of our Theorems A and B it will follow that $Z_G$ is a cone.
If one is able to prove this directly, the proof of Theorem A follows easily and directly 
from Lemma 6 (see Theorem 8). It is very easy to show that $Z_G$ is a cone if and only if $th$ is in fact
the full polar map $ \partial G$ if and only if in $ \Bbb P^{n *}$ the point corresponding 
to $H_0$ is not on a line connecting $Z_{L_i}^{*}$ and $Z_{L_j}^{*}$ for $i \neq j \neq 0$ 
if and only if the polynomial $G_0$ is square free.
This is the hard part of the problem, connected with cohomological properties of the corresponding arrangement of hyperplanes in $\Bbb P^n$ and with syzygies of the base locus of $ \partial F$. It is in fact remarkable that Dolgachev proves it directly in case $n=3$.
\vskip 0.2 cm
Let us now consider homaloidal polynomials which are products of linear forms with at least a square factor; notice that such polynomials arise naturally from square free homaloidal polynomials which are products of linear forms: in Lemma 6 we have proven that the restriction of a homaloidal square free product of linear forms induces on each component
of $Z_F$ a homaloidal product of linear forms.
Quite surprisingly {\it there is  a priori no relation among \/} $ \partial F$ and
$ \partial F_{red}$.
\vskip 0.2 cm
Let us choose in fact $r+1$ distinct linear forms $L_0,\ldots,L_r$ in $\Bbb P^n$ together with an identification $X_0=L_0$, and consider the following polynomials, where $H_0$ is the hyperplane of equation $X_0=0$ and $L_{i,0}=L_i \cap H_0$:
$$F=X_0^{m_0}\prod_{i=1}^r L_i^{m_i}, \qquad F^{'}=X_0^{m_0-1}\prod_{i=1}^r L_i^{m_i-1},
\qquad F_{red}= \frac {F}{F^{'}}, $$
$$G=\prod_{i=1}^r L_i^{m_i}, \qquad G^{'}=\prod_{i=1}^r L_i^{m_i-1}, \qquad G_{red}= \frac {G}{G^{'}},$$
$$G_0=G \cap H_0, \qquad  G_0^{'}= G^{'} \cap H_0.$$
We compute the moving parts of the polar systems defined by $F$ and $F_{red}$. We have:
$$ | M( \partial F)|=<m_0 G_{red}+X_0 \sum_{i=1}^r m_i \frac {\partial L_i}
{\partial X_0} \prod_{j \neq i,0} L_j, \ldots, X_0 \sum_{i=0}^r m_i \frac {\partial L_i}
{\partial X_n} \prod_{j \neq i,0} L_j >|,$$
$$ | M( \partial F_{red})|= | \partial F_{red}|= <G_{red}+X_0 \sum_{i=1}^r  \frac {\partial L_i}
{\partial X_0} \prod_{j \neq i,0} L_j, \ldots, X_0 \sum_{i=0}^r  \frac {\partial L_i}
{\partial X_n} \prod_{j \neq i,0} L_j >|.$$
We also compute an explicit basis of a system $|G_0^{''}|$ which sits in a chain $|M( \partial G_0| \subset |G_0^{''}| \subset | \partial G_0|$ on $H_0$:
$$ |G_0^{''}|=< \sum_{i=1}^r m_i \frac {\partial L_i}
{\partial X_1} \prod_{j \neq i} L_{j,0}, \ldots,\sum_{i=1}^r m_i \frac {\partial L_i}
{\partial X_n} \prod_{j \neq i} L_{j,0}>|,$$

Consider now the following diagram of maps, where $f$ and $g$ induce a minimal resolution of the morphism induced by $|M( \partial F)|$:

$$
\CD
  \,    X       \\
  \llap {$f$} \swarrow  \qquad \qquad \searrow \rlap {$g$} \\
\, \qquad H_0 \subset \Bbb P^n --{M (\partial F)}--\to    \Bbb P^{ n *}-\pi - \to  P. \\
\endCD
$$

We define $D$ to be the strict transform of $H_0$ in $X$, we denote by $Y$ a general member of $|M( \partial F)|$, by $Y^{'}$ a general member of the system $|M ( \partial G)|$, by
$\overline {Y}$ and $\overline {Y^{'}}$ their strict transforms on $X$.
Quite surprisingly, all the arguments used in order to prove Lemma 5 and Lemma 6 apply verbatim and in particular it holds the following: 
\proclaim {Lemma 7}
With the above notations:
\roster
\item $D$ is the strict transform on $X$ of the exceptional divisor in the blow up  $h^{'} : Z \to \Bbb P^{n*}$ of $\Bbb P^{n*}$ at $U_0$,
\item the restriction of the linear system $|(th)^{*} \Cal O_P(1)|$  to $D$ induces a morphism 
$m: D \to P$ which is a resolution of singularities of the polar map defined by $G_0$ on 
$H_0$, i.e. $|(\overline { Y} -D)|_D| = |M ( \partial G_0)|$, 
\item the polynomial $G_0$ is homaloidal in $H_0$,
\endroster
\endproclaim
\demo {Proof} \qed\enddemo

We are now able to prove Theorems A and B at once.
\proclaim { Theorem 8}
Let $ L_0, \ldots, L_r$ be distinct linear forms and let $F=\prod_{i=0}^r L^{m_i}$, with $m_i \geq 1$ for all $i=0, \ldots, r$.
Then $F$ is homaloidal if and only if $F_{red}=\prod_{i=0}^r L_i$ is homaloidal
and $F_{red}$ is homaloidal if and only if $r=n$ and the $L_i 's$ are independent linear forms .
\endproclaim
\demo {Proof}
The proof is by induction on $n$.

The starting point of the induction is the case $n=1$ which is easy: if
$F=\prod_{i=0}^r L_i^{m_i}$ is homaloidal the base point free system $| M ( \partial F)|$ must be of degree one, from which it follows easily that $r=1$ and that $L_0$ and $L_1$ are
in linear general position ( they are distinct by hypothesis). The converse is, up to a projectivity, been proved after Definition 3. The same argument works a fortiori if $F$ is square free.

Let us then move to $\Bbb P^n$, with $ n > 1$ and consider first the case of a square free homaloidal polynomial $F=\prod_{i=0}^r L_i$. We plug $X_0=L_0$, and we apply Lemma 6 in order to get a homaloidal polynomial $G_0=\prod_{i=1}^r L_{i,0}$ on $H_0$. If $G_0$ is reduced we have that
by induction $r=n$ and the $L_{i,0} 's$ are independent in $H_0$, from which it follows that  $X_0,L_1, \ldots, L_n$ are independent in $ \Bbb P^n$. Consider then the case in which $G_0$ is not reduced. $G_0$ is still homaloidal so that by induction we can reorder the $L_i 's$ in such a way that $L_{1,0}, \ldots, L_{n,0}$ are independent and there exist $m_1, \ldots, m_n$ for which $ \sum_{i=1}^n m_i = r$ and
$$G_0= \prod_{i=1}^n L_{i,0}^{m_i},$$
in such a way that
\roster
\item $L_{n+1}, \ldots, L_{n+m_1-1}$ are in $<X_0,L_1>$,
\item $\cdot \cdot$
\item $L_{n+m_1+ \ldots m_{n-1}-n+2}, \ldots, L_{n+\sum_{i=1}^nm_i-n}$ are in $<X_0,L_n>$.
\endroster
In other words, looking at the dual points in $\Bbb P^{n*}$, we must have that all $Z_{L_i}^{*}$'s, with $i \geq n+1$ must lie in the cone with vertex $U_0$ projecting $Z_{L_1}^{*}, \ldots, Z_{L_n}^{*}$. But we can apply the same reasoning we have applied to $L_0=X_0$ to any other $L_i=X_i$ for all $i=1,\ldots,n$ and the intersection of all these cones is empty, so that there will exists some $i \in \{0,\ldots,n \}$ for which the corresponding homaloidal polynomial $G_i$ is reduced. We then proceed as if $G_0$ were reduced.

Suppose now that $F=\prod_{i=0}^r L_i^{m_i}$ is nonreduced and homaloidal. We must prove that
$ F_{red}$ is homaloidal, the converse being a consequence of the first part of this Proof and of the example after Definition 3.
Plugging $X_0=L_0$ and applying Lemma 7 we get that $| M ( \partial F)|$ induces on $H_0$ the homaloidal system defined by $ G_0$. By induction and by the same argument as above, we get the thesis, i.e. that
 $r=n$ and the linear forms $L_0, \ldots,L_n$ are independent, so that
$ F_{red} $ is homaloidal.\enddemo\qed

Andrea Bruno, Universita' degli Studi Roma Tre, Dipartimento di Matematica, Largo S. L: Murialdo, 1, 00146 Roma (Italia)\newline
E-mail address: bruno at mat.uniroma3.it

\end
\widestnumber \key{DiPa}
\Refs \nofrills {References}

\ref \key Di
\by Dimca,A.
\pages 47--53
\paper On polar Cremona transformations
\yr 2001 \vol n.1
\jour An. Stiinct. Univ. Ovidius Constancta Ser. Mat. 9 
 \endref


\ref \key DiPa
\by Dimca,A.; Papadima,S.
\pages 473-507
\paper Hypersurface complements, Milnor fibers and higher homotopy groups of arrangments
\yr 2003 \vol no. 2  vol 158
\jour Ann. of Math. (2)
 \endref


\ref \key Do
\by Dolgachev,I.
\pages 191--202
\paper Polar Cremona transformations
\yr 2000 \vol 48 
\jour Michigan Math. J.(volume dedicated to William Fulton)
 \endref
 
 \ref \key EKP
\by Etingof,P.;Kazhdan,D.;Polishchuk, A.
\pages 27--66
\paper When is the Fourier transform of an elementary function elementary?
\yr 2002 \vol 8, n.1
\jour Selecta Math (N.S.)
 \endref
 
\ref \key KiSa
\by Kimura,T.;Sato,M.
\pages 1-155
\paper A classification of prehomogeneous vector spaces and their relative invariants
\yr 1977 \vol 65
\jour Nagoya Math. Journal
 \endref
 
\ref \key KS
\by Kraft,H;Schwarz,G.
\pages 641-649
\paper Rational covariants of reductive groups and homaloidal polynomials
\yr 2001 \vol 8
\jour Mathematical Research Letters
\endref

\ref \key R
\by Russo,F.
\paper Tangents and Secants of algebraic varieties. Notes of a course 24th Coloquio Brasileiro de Matematica, IMPA
\yr 2003
\endref



\endRefs
